\documentclass[12pt]{article}
\usepackage{amsxtra,amssymb,amsthm,amsmath,latexsym}

\theoremstyle{plain}

\newtheorem{theorem}{Theorem}

\numberwithin{equation}{section}

\def\oH{\buildrel\circ\over H}
\def\oH1{\buildrel\circ\over H\kern-.02in{}^1}
\def\qed{{\hfill $\Box$}}

\begin{document}

\title{Finding discontinuities of piecewise-smooth functions}

\author{
A.G. Ramm\\
 Mathematics Department, Kansas State University, \\
 Manhattan, KS 66506-2602, USA\\
ramm@math.ksu.edu}

\date{}

\maketitle\thispagestyle{empty}

\begin{abstract} \footnote{Math subject classification:  
65D35; 65D05} \footnote{key words: stable differentiation,
noisy data, discontinuities, jumps, signal processing, edge
detection.}

Formulas for stable differentiation of piecewise-smooth
functions are given. The data are noisy values of these
functions. The locations of discontinuity points and the
sizes of the jumps across these points are not assumed
known, but found stably from the noisy data. 

\end{abstract}


\section{Introduction}
Let $f$ be a piecewise-$C^2([0,1])$ function,
$0<x_1<x_2<\dots <x_J, 1\leq j\leq J$,
are discontinuity points of $f$. We do not assume their
locations $x_j$ and their number $J$ known a priori.
We assume that the limits $f(x_j\pm 0)$ exist,
and
\begin{equation}\label{e1.1}
 \sup_{x\not= x_j, 1\leq j\leq J} |f^{(m)} (x)| \leq M_{m},
 \quad m=0,1,2.
 \end{equation}
Assume that $f_\delta$ is given,
$\|f-f_\delta\|:=\sup_{x\not= x_j, 1\leq j\leq J} |f-f_\delta|\leq 
\delta$,
where $f_\delta\in L^\infty(0,1)$ are the noisy data.

{\it The problem is: given $\{f_\delta,\delta\},$
where $\delta\in (0,\delta_0)$ and $\delta_0>0$ is a 
small number, 
estimate
stably $f'$, find the locations of
discontinuity points  $x_j$ of $f$ and their number $J$,
and estimate the jumps $p_j:=f(x_j+0)-f(x_j-0)$ of $f$ 
across $x_j$,
$1\leq j \leq J$.}

A stable estimate $R_\delta f_\delta$ of $f'$ is an estimate 
satisfying the relation $\lim_{\delta \to 0}||R_\delta 
f_\delta-f'||=0$.

There is a large literature on stable differentiation of noisy
smooth functions
(e.g., see references in \cite{R470}),
but the problem stated above was not solved for 
piecewise-smooth functions by the method
given below. A statistical estimation of the 
location of discontinuity
points from noisy discrete data is given in \cite{R316}.

The following formula was proposed originally (in 1968, see 
references in \cite{R470}) for
stable estimation of $f'(x)$, assuming
$f\in C^2([0,1])$,  $M_2\neq 0$, and given noisy data $f_\delta$:
\begin{equation}\label{e1.2}
 R_\delta f_\delta:=
 \frac{f_\delta (x+h(\delta))-f_\delta(x-h(\delta))}{2h(\delta)},\,\,
 h(\delta):=
 \left( \frac{2\delta}{M_2}\right)^{\frac{1}{2}},
 \quad
 h(\delta)\leq x\leq 1-h(\delta),
 \end{equation}
 and
\begin{equation}\label{e1.3}
 \|R_\delta f_\delta -f'\| \leq \sqrt{2M_2\delta}:=\varepsilon(\delta),
 \end{equation}
where the norm in (1.3) is 
$L^\infty(0,1)-$norm.
Moreover, (cf \cite{R470}),
\begin{equation}\label{e1.4}
 \inf_T \sup_{f\in K(M_2,\delta)}
 \|Tf_\delta-f'\|
 \geq \varepsilon(\delta),
 \end{equation}
where $T:L^\infty(0,1)\to L^\infty(0,1)$ runs through the
set of all bounded operators,
$K(M_2,\delta):=\{f:\|f''\|\leq M_2,\ \|f-f_\delta\|\leq\delta\}$.
Therefore estimate \eqref{e1.2} is the best possible estimate
of $f'$, given noisy data $f_\delta$,
and assuming $f\in K(M_2,\delta)$.

In \cite{R470} this result was generalized to the case $f\in
K(M_a,\delta)$, $\|f^{(a)}\|\leq M_a$, $1<a\leq 2$, where
$\|f^{(a)}\|:=\|f\|+ \|f'\| +\sup_{x,x'} \frac{
|f'(x)-f'(x')| } { |x-x'|^{a-1} }$, $1<a\leq 2$, and
$f^{(a)}$ is the fractional-order derivative of $f$.

The aim of this paper is to extend the above results to the
case of piecewise-smooth functions. In Section 2 the results
are formulated, and  proofs are given. 
In Section 3
the case of continuous piecewise-smooth functions is
treated.

\section{Formulation of the result}

\begin{theorem}\label{T:1}
Formula \eqref{e1.2} gives stable estimate of $f'$ on the set
$S_\delta:=[h(\delta), 1-h(\delta)]
\,\backslash\, \bigcup_{j=1}^J (x_j-h(\delta), x_j+h(\delta))$,
and \eqref{e1.3} holds with the norm $\|\cdot\|$
taken on the set $S_\delta$. Assuming $M_2>0$ and computing the quantities 
$f_j:=\frac{f_\delta (jh+h)-f_\delta(jh-h)}{2h}$, where $ h:=h(\delta):=
 \left( \frac{2\delta}{M_2}\right)^{\frac{1}{2}}$, 
$1\leq j <[\frac 1 h],$
for sufficiently small $\delta,$ one finds the location of discontinuity
points of $f$ with accuracy $2h$, and their number $J$. Here 
$[\frac 1 h]$ is the integer smaller than $\frac 1 h$ and closest to 
$\frac 1 h$. The discontinuity points of $f$ are located on the intervals 
$(jh-h, jh+h)$ such that $|f_j|\gg 1$ for sufficiently 
small $\delta$, where $\varepsilon (\delta)$ is defined in (1.3).
The size $p_j$ of the jump of $f$ across the discontinuity point $x_j$
is estimated by the formula $p_j\approx f_\delta 
(jh+h)-f_\delta(jh-h)$, and the error of this estimate is 
$O(\sqrt 
{\delta})$. 
\end{theorem}

Let us assume that $\min_j p_j:=p\gg h(\delta)$, where
$\gg$ means "much greater than".
Then $x_j$ is located on the $j\hbox{-th}$ interval
$[jh-h,jh+h]$, $h:=h(\delta)$, such that
\begin{equation}\label{e2.1}
|f_j|:= \left| \frac{f_\delta(jh+h) - f_\delta(jh-h)}
  {2h} \right|
 \gg 1,
 \end{equation}
so that $x_j$ is localized with the accuracy $2h(\delta)$.
More precisely, 
$|f_j|\geq \frac {|f(jh+h)-f(jh-h)|}{2h}-\frac {\delta}{h}$,
and $\frac {\delta}{h}=0.5\varepsilon(\delta)$, where 
$\varepsilon(\delta)$ is defined in (1.3).  One has
$$|f(jh+h)-f(jh-h)|\geq 
p_j-|f(jh+h)-f(x_j+0)|-|f(jh-h)-f(x_j-0)|\geq
p_j-2M_1h.$$
 Thus,
$$|f_j|\geq \frac 
{p_j}{2h}-M_1-0.5\varepsilon(\delta)=c_1\frac{p_j}{\sqrt{\delta}}-c_2\gg 
1,$$
where $c_1:=\frac{\sqrt{M_2}}{2\sqrt{2}}$, and 
$c_2:=M_1+0.5\varepsilon(\delta)$.
 
The jump
$p_j$ is estimated by the formula: 
\begin{equation}\label{e2.2}
p_j\approx [f_\delta(jh+h) - f_\delta(jh-h)],
\end{equation}
 and the error estimate of this formula can be given:
\begin{equation}\label{e2.3}
 |p_j-[f_\delta (jh+h)-f_\delta(jh-h)]|
 \leq 2\delta+2M_1h=2\delta+2M_1
 \sqrt{\frac{2\delta}{M_2}}=O(\sqrt{\delta}).
 \end{equation}
Thus, the error of the calculation of $p_j$ by the formula 
$p_j\approx f_\delta (jh+h)-f_\delta(jh-h)$ is 
  $O(\delta^{\frac{1}{2}})$
as $\delta\to 0$.

{\bf Proof of Theorem 1.}
If $x\in S_\delta$, then using Taylor's formula one gets:
\begin{equation}\label{e2.4}
 |(R_\delta f_\delta)(x)-f'(x)|\leq \frac{\delta}{h}+\frac{M_2h}{2}.
 \end{equation}
Here we assume that $M_2>0$ and the interval 
$(x-h(\delta),x+h(\delta))\subset S_\delta,$ i.e., this interval
does not contain discontinuity points of $f$.
If, for all sufficiently small $h$, not necessarily for $h=h(\delta)$,
inequality (2.4) fails, i.e., if $|(R_\delta f_\delta)(x)-f'(x)|>
\frac{\delta}{h}+\frac{M_2h}{2}$ for all sufficiently small $h>0$,
then the interval $(x-h,x+h)$ contains a point $x_j\not\in 
 S_\delta,$ i.e., a point of discontinuity of $f$ or $f'$. This 
observation can be used for locating the position of an
isolated  discontinuity point $x_j$ of $f$ with any desired accuracy
provided that the size $|p_j|$ of the jump of $f$ across $x_j$
is greater than $4\delta$,  $|p_j|> 4\delta$, and that $h$ can be taken as 
small as desirable. Indeed, if  $x_j\in (x-h,x+h)$, then we have 
$2\delta +M_2h^2<|f_\delta(x+h)-f_\delta(x-h)-2hf'(x)|=
|p_j+O(h^2)\pm 2\delta|$. The above estimate follows from the relation
$|f_\delta(x+h)-f_\delta(x-h)-2hf'(x)|=|f(x+h)-f(x_j+0)+p_j+
f(x_j-0)-f(x-h)\pm 2\delta|=|p_j+O(h^2)\pm 2\delta|$.
Here $|p\pm b|$, where $b>0$, denotes  a quantity such that $|p|-b\leq
|p\pm b|\leq |p|+b$.
 Thus, if $h$ is sufficiently small and  $|p_j|> 4\delta$,
then the inequality $2\delta +M_2h^2<|p_j+O(h^2)\pm 2\delta|$ can be 
checked, and 
therefore the inclusion $x_j\in (x-h,x+h)$ can be checked. Since $h>0$
is arbitrarily small in this argument, it follows that the location
of the discontinuity point $x_j$ of $f$ is established with arbitrary
accuracy.
Additional discussion of the case
when a discontinuity point $x_j$ belongs to the interval 
 $(x-h(\delta),x+h(\delta))$ will be  given below. 

Minimizing the right-hand side of \eqref{e2.4} with respect to
$h$ yields formula \eqref{e1.2}  for the 
minimizer $h=h(\delta)$ defined
in (1.2), and estimate \eqref{e1.3} for the minimum of the 
right-hand side of \eqref{e2.4}.

If $p\gg h(\delta)$, and \eqref{e2.1} holds,
then the discontinuity points are located with the accuracy
$2h(\delta)$, as we prove now. 

Consider the case when a discontinuity point $x_j$ of $f$ belongs
to the interval $(jh-h,jh+h)$, where $h=h(\delta)$. Then
estimate \eqref{e2.2} can be obtained as follows.
For $jh-h\leq x_j\leq jh+h$, one has
$$\begin{aligned}
 |f(x_j+0)
 &-f(x_j-0) -f_\delta(jh+h) +f_\delta(jh-h)| \leq 2\delta+ \\
 & +|f(x_j+0)-f(jh+h)| \\
 & \qquad + |f(x_j-0)-f(jh-h)| \leq 2\delta + 
M_12h,\quad h=  h(\delta).
 \end{aligned}
$$ 
This yields formulas (2.2) and (2.3). Computing the
quantities $f_j$ for $1\leq j <[\frac 1 h]$, and finding the
intervals on which (2.1) holds for sufficiently small
$\delta$, one finds the location of discontinuity points of
$f$ with accuracy $2h$, and the number $J$ of these points.  
For a small fixed $\delta>0$ the above method allows one to recover
the discontinuity points of $f$ at which $|f_j|\geq \frac 
{|p_j|}{2h}- \frac{\delta}{h}-M_1\gg 1$. This is the inequality (2.1).
 If $h=h(\delta)$, 
then $\frac{\delta}{h}=0.5 \varepsilon(\delta)= O(\sqrt{\delta})$, and 
 $|2hf_j-p_j|= O(\sqrt{\delta})$ as $\delta \to 0$ provided that $M_2>0$.
Theorem 1 is proved. \qed

{\bf Remark 1:}
Similar results can be derived if
$\|f^{(a)}\|_{L^{\infty}(S_\delta)}:=\|f^{(a)}\|_{S_\delta} \leq M_a$, 
$1<a\leq 2$.
In this case
$ h=h(\delta)=c_a \delta^{\frac{1}{a}}$, where
$c_a=\left[ \frac{2}{M_a(a-1)} \right]^{\frac{1}{a}}$,
$R_\delta f_\delta$ is defined in \eqref{e1.2}, and the 
error 
of the estimate is:
$$
 \|R_\delta f_\delta -f'\|_{S_\delta}
 \leq a M^{\frac{1}{a}}_a
 \left( \frac{2}{a-1} \right)^{\frac{a-1}{a}} \delta^{\frac{a-1}{a}}.
 $$
The proof is similar to that given in Section 3.
It is proved in \cite{R470} that for $C^a\hbox{-functions}$ given
with noise it is possible to construct stable differentiation
formulas if $a>1$ and it is impossible to construct
such formulas if $a\leq 1$.
The obtained formulas are useful in applications.
One can also use $L^p\hbox{-norm}$ on $S_\delta$ in the estimate
$\| f^{(a)}\|_{S_\delta} \leq M_a$ (cf. \cite{R470}).

{\bf Remark 2:} The case when $M_2=0$ requires a special discussion. In 
this case 
the last term on the right-hand side of formula (2.4) vanishes
and the minimization with respect to $h$ becomes void: it requires that 
$h$ be as large as possible, but one cannot take $h$ arbitrarily large 
because estimate (2.4) is valid only on the interval $(x-h,x+h)$ which 
does not contain discontinuity points of $f$, and these points are
unknown. 
If $M_2=0$, then $f$ is a piecewise-linear function.  
The discontinuity points of a piecewise-linear function 
can be found if the sizes $|p_j|$ of the jumps of $f$ across these points 
satisfy the inequality
$|p_j|>>2\delta +4M_1h$ for some choice of $h$. For instance, 
if $h=\frac{\delta}{2M_1}$, then $2\delta +4M_1h=4\delta.$
So, if $|p_j|>>4\delta,$ then the location of discontinuity points of $f$ 
can be found  in the case when $M_2=0$.
 These  points are located on the intervals for which 
$|f_\delta(jh+h) - f_\delta(jh-h)|>>4\delta$, where 
 $h=\frac{\delta}{2M_1}$.

The size $|p_j|$ of the jump of $f$ across a discontinuity point $x_j$
can be estimated by  formula (2.2) with $h=\frac{\delta}{2M_1}$, and 
one assumes that $x_j\in (jh-h,jh+h)$ is the only discontinuity point on 
this interval.
The  error of the formula (2.2) is estimated as in the proof of Theorem 1.
This error is not more than $2\delta+4M_1h=4\delta$ for the above choice 
of $h=\frac{\delta}{2M_1}$.

One can estimate the derivative of $f$ at the point of smoothness of $f$ 
assuming  $M_2=0$ provided that this derivative is 
not too small. If $M_2=0$, then $f=a_jx+b_j$ on every interval $\Delta_j$
between the discontinuity points $x_j$, where $a_j$ and $b_j$ are some 
constants. If $(jh-h,jh+h)\subset \Delta_j$, and 
$f_j:=\frac{f_\delta(jh+h)-f_\delta(jh-h)}{2h}$, then
$|f_j-a_j|\leq \frac {\delta}{h}$. Choose $h=\frac {t\delta} {M_1},$
 where $t>0$ is a parameter, and $M_1=\max_{j}|a_j|.$ Then 
the relative error of the approximate formula  $a_j\approx f_j$
for the derivative $f'=a_j$
on $\Delta_j$  equals to $\frac {|f_j-a_j|}{|a_j|}\leq 
\frac {M_1}{t|a_j|}$. Thus, if, e.g., $|a_j|\geq \frac {M_1}2$ and 
$t=10$,
then the relative error of the above approximate formula is not more than 
$0.1$.
  
\section{Continuous piecewise-smooth functions} 
Suppose now
that $\xi \in (mh-h, mh+h)$, where $m>0$ is an integer, and
{\it $\xi$ is a point at which $f$ is continuous but
$f'(\xi)$ does not exist}. Thus, the jump of $f$ across
$\xi$ is zero, but $\xi$ is not a point of smoothness of
$f$.  {\it How does one locate the 
point $\xi$?}

The algorithm we propose consists of the following.
 We assume that $M_2>0$ on $S_\delta$.
Calculate the numbers
$f_j:=\frac{f_\delta(jh+h)-f_\delta(jh-h)}{2h}$ and $|f_{j+1}-f_j|$,
$j=1,2,\dots$, $h=h(\delta)=\sqrt{\frac{2\delta}{M_2}}$.
Inequality (1.3) implies $f_j-\varepsilon(\delta)\leq f'(jh)\leq 
f_j+\varepsilon(\delta)$, where $\varepsilon(\delta)$ is defined
in (1.3). 

Therefore, if $ |f_j|>\varepsilon(\delta)$,
then $\hbox{sign } f_j=\hbox { sign }f'(jh)$. 

One has:
$$J-\frac {\delta}{h}\leq |f_{j+1}-f_j|\leq J+\frac {\delta}{h},$$
 where 
$\frac {\delta}{h}=0.5 \varepsilon(\delta)$ and $J:=|\frac 
{f(jh+2h)-f(jh)-f(jh+h)+f(jh-h)}{2h}|.$ 
Using Taylor's formula, one derives the estimate:
\begin{equation}\label{e3.1} 0.5 [J_1-\varepsilon(\delta)]\leq J\leq  0.5 
[J_1+\varepsilon(\delta)],
\end{equation}
where $J_1:= |f'(jh+h)-f'(jh)|$.

If the interval $(jh-h,jh+2h)$ belongs to $S_\delta$, then
$J_1=|f'(jh+h)-f'(jh)|\leq M_2h=\varepsilon(\delta)$. In this case
$J\leq \varepsilon(\delta)$, so 
\begin{equation}\label{e3.2}|f_{j+1}-f_j|\leq \frac 3 2 
\varepsilon(\delta)
\quad \hbox { if } \quad (jh-h,jh+2h)\subset S_\delta.
\end{equation}

{\it Conclusion: if $|f_{j+1}-f_j|> \frac 3 2 
\varepsilon(\delta)$,
then the interval $(jh-h,jh+2h)$ does not belong to $S_\delta$,
that is, there is a point $\xi\in (jh-h,jh+2h)$ at which 
the function $f$ is not twice continuously differentiable with 
$|f^{\prime \prime}|\leq M_2$. Since 
we assume that 
either at a point $\xi$ the function is twice differentiable, or at this 
point $f'$ does not exist, it follows that if
$|f_{j+1}-f_j|>\frac 3 2 \varepsilon(\delta)$, then there is 
a point 
$\xi\in (jh-h, jh +2h)$ at which $f'$ does not exist.} 

If  
\begin{equation}\label{e3.3}
 f_jf_{j+1}<0,
 \end{equation}
and
\begin{equation}\label{e3.4}
 \min (|f_{j+1}|,|f_j|)>\varepsilon(\delta),
 \end{equation}
then (3.3) implies $f'(jh)f'(jh+h)<0$, so the interval $(jh,jh+h)$ 
contains a critical point $\xi$ of $f$,
or a point $\xi$ at which $f'$ does not exist. To determine which
one of these two cases holds, let us use the right inequality (3.1).
If $\xi$ is a critical point of $f$ and $\xi\in (jh, jh+h)\subset 
S_\delta$,  then $J_1\leq 
\varepsilon(\delta)$, and in this case the right inequality (3.1) 
yields
\begin{equation}\label{e3.5} 
 |f_{j+1}-f_j|\leq  \varepsilon(\delta).  
 \end{equation}
{\it Conclusion: If (3.3)-(3.5) hold, then $\xi$ is a critical point.
If (3.3) and (3.4) hold and $ |f_{j+1}-f_j|>  
\varepsilon(\delta)$ then  $\xi$ is a
point of discontinuity of $f'$.}

If $\xi$ is a point of discontinuity of $f'$, we would like to estimate 
the jump 
$$P:=|f'(\xi+0)-f'(\xi-0)|.$$
Using Taylor's formula one gets 
\begin{equation}\label{e3.6}
f_{j+1}-f_j=\frac P{2}\pm 3.5 \varepsilon(\delta). 
 \end{equation}
The expression $A=B\pm b,\,\, b>0,$ means that $B-b\leq 
A\leq B+b$.
Therefore, 
\begin{equation}\label{e3.7}
P=2(f_{j+1}-f_j) \pm 7 \varepsilon(\delta).
 \end{equation}
We have proved the following theorem:

\begin{theorem}\label{T:2}
If $\xi\in (jh-h,jh+2h)$ is a point of continuity of $f$ and 
$|f_{j+1}-f_j|>\frac 3 2 \varepsilon(\delta)$, then $\xi$ is a point 
of 
discontinuity of $f'$. If
 \eqref{e3.3} and \eqref{e3.4} hold, and 
$|f_{j+1}-f_j|\leq \varepsilon(\delta)$, then $\xi$ is a critical point
of $f$. If  \eqref{e3.3} and \eqref{e3.4} hold and $|f_{j+1}-f_j|> 
\varepsilon(\delta)$, then $\xi\in (jh,jh+h)$ is a point of 
discontinuity of $f'$. The jump $P$ of $f'$ across $\xi$ is 
estimated by formula (3.7). 
\end{theorem}

\end{document}